\documentstyle{article}
\begin{document}

\newcommand{\qq}{\mbox{$q^{\prime}$}}
\newcommand{\pp}{\mbox{$p^{\prime}$}}
\newcommand{\pip}{\mbox{${\scriptstyle p}'$}}
\newcommand{\qiq}{\mbox{${\scriptstyle q}'$}}
\newcommand{\ff}{\mbox{$f^{\prime}$}}
\newcommand{\DD}{\mbox{$D^{\prime}$}}
\newcommand{\DeDe}{\mbox{$\Delta^{\prime}$}}
\newcommand{\tautau}{\mbox{$\tau^{\prime}$}}
\font\tenBbb=msbm10  \font\sevenBbb=msbm7  \font\fiveBbb=msbm5
\newfam\Bbbfam
\textfont\Bbbfam=\tenBbb \scriptfont\Bbbfam=\sevenBbb
\scriptscriptfont\Bbbfam=\fiveBbb
\def\Bbb{\fam\Bbbfam\tenBbb}
\def\R{{\Bbb R}}
\def\C{{\Bbb C}}
\def\N{{\Bbb N}}
\def\zed{{\Bbb Z}}

\title{A UNIFORMLY CONVEX HEREDITARILY INDECOMPOSABLE BANACH SPACE}
\author{Valentin Ferenczi}
\date{October 1994 }

\maketitle

\begin{abstract}
We construct a uniformly convex hereditarily indecomposable Banach 
space, 
  using similar methods as Gowers and Maurey in \cite{GM}
 and the theory
 of complex
interpolation for a family of Banach spaces of Coifman, Cwikel,
 Rochberg, Sagher and
Weiss ( \cite{CC}).
\end{abstract}

\paragraph{Introduction}

A {\em hereditarily indecomposable (or H.I.)} space is an
infinite dimensional Banach space such that no
 subspace can
be written as the topological sum of two infinite dimensional
subspaces. As an easy consequence,
 no such space can contain an unconditional basic sequence. This
notion also appears as the 'worst' type of subspace
 of a Banach space
 in \cite {G}.
In \cite{GM}, Gowers and Maurey constructed the first
known example of a hereditarily indecomposable space.
Gowers-Maurey space is reflexive, however it is not
uniformly convex. In this article, we provide an example of
a uniformly convex hereditarily indecomposable space.

 \section{A class of uniformly convex Banach spaces}
\subsection{Definitions}

Let $c_{00}$ be the space of sequences of scalars all but finitely many of
which are zero. Let $e_{1},e_{2},\ldots$ be its unit vector basis.
If $E\subset \N$, then we shall also use the letter
$E$ for the projection from $c_{00}$ to $c_{00}$ defined by
$E(\sum_{i=1}^\infty a_i e_i)=\sum_{i\in E}a_i e_i$. If
$E,F \subset \N$, then we write $E<F$ 
to mean that $\sup E<\inf F$.
 An {\em interval} of integers is a subset of
$\N$ of the form $\{a,a+1,\dots,b\}$ for some $a,b\in \N$.
For $N$ in $\N$, $E_{N}$ denotes the interval
$\{1,\ldots,N\}$. The {\em range} of a vector
 $x$ in $c_{00}$, written $ran(x)$, is the smallest interval $E$
such that $Ex=x$.
  We shall write $x<y$ to mean $ran(x)<ran(y)$. If
$x_1<\cdots<x_n$ we shall say that $x_1,\ldots,x_n$ are
 {\em successive}.

The corresponding
 notation about range and successive functions
 will be used for analytic functions with values in $c_{00}$
( the range of such functions is always finite).
 Let $\cal X$ be the class of
 normed spaces of the form $(c_{00},\|.\|)$, such
that $(e_{i})_{i=1}^{\infty}$ is a normalized bimonotone basis.
By a {\em block basis} in a space $X \in {\cal X}$ we mean a sequence
$x_1,x_2,\dots$ of successive non-zero vectors in $X$  ( note that such a
sequence must be a basic sequence) and by a {\em block subspace}
 of a space
$X\in\cal X$ we mean a subspace generated by a block basis.

Let $f$ be the function $\log_{2}(x+1)$.
If 
 $X \in {\cal X}$, and all successive vectors $x_{1}, \ldots,
x_{n}$ in $X$ satisfy 
the inequality $f(n)^{-1} \sum_{i=1}^{n} \|x_{i}\| \leq
 \| \sum_{i=1}^{n} x_{i} \|$,
then we say that $X$ satisfies a {\em lower $f$-estimate}.

Let $q >1$ in $\R$, $\qq$ such that 
$1/q + 1/\qq=1$. Let $\theta \in \ ]0,1[\ $, and 
$p$ be the number defined by $ 1/p = 1-\theta + \theta /q $.

Let $S$ be the strip $\{ z \in \C / Re(z) \in [0,1] \}$,
 $\delta S$ its boundary, $S_{0}$
the line $\{z / Re(z)=0 \}$, $S_{1}$ the line $\{ z/Re(z)=1 \}$.
 Let $\mu$ be the Poisson probability measure associated to the point $\theta$
for the strip $S$. We have $\mu(S_{0})=1-\theta$.
Let $\mu_{0}$ be the probability measure on $\R$
 defined by $\mu_0(A)=\mu(iA)/(1-\theta)$,
 $\mu_{1}$ be the probability measure on $\R$ defined by
 $\mu_1(A)=\mu(1+iA)/\theta$.
Let ${\cal A}_{S}$ be the set of analytic functions $F$
on $S$, with values in $c_{00}$, which
are $L_{1}$ on $\delta S$ for
$d\mu$ and which satisfy the
 Poisson integral
representation $F(z_0)=\int_{\delta S}F(z)dP_{z_0}(z)$
 on $S$
( this is well defined since such functions
 have finite ranges). If $F$ is analytic and 
bounded on $S$, then $F \in {\cal A}_S$.

We recall the definition of the interpo\-lation space of a
 fa\-mily of
$N$-dimen\-sional spaces from \cite{CC}. Let
$\|.\|_z$ for $z$ in $S$ be a family of norms
 on $\C^N$, equivalent 
with log-integrable constants, and such that $z \mapsto \|x\|_{z}$ is
 measurable for all $x$ in $\C^N$.
 The interpolation space in $\theta$ is
defined by the norm 
 $ \|x\| = \inf_{F \in {\cal A}_S^N,\ F(\theta)=x}
 (\int_{z \in \delta S} \|F(z)\|_{z}  d \mu (z))$,
where ${\cal A}_S^N$ denotes the image of
the canonical projection
from ${\cal A}_S$ into the space of functions from $S$ to
$\C^N$.

We generalize to the infinite-dimensional case as follows.
Let $\{ X_z, z \in \delta S \}$ be
 a family of Banach spaces in $\cal X$, equipped
with norms $\|.\|_z$, such that 
for all $x$ in $c_{00}$, the function $z \mapsto \|x\|_{z}$ is
 measurable, and such that over vectors of finite range $N$,
 the norms $\|.\|_z$ are equivalent with log-integrable
constants.
Let
$X_z^N$ be $E_N X_z$, $X^{N}$ be the $\theta$-interpolation
space
 of 
the family $X_z^N$; the interpolation space of the family in 
$\theta$ is
 $ completion( \cup_{N \in \N} X^{N})$.

Now let $\{X_{t}, t \in \R \}$ be a family
 of spaces in $\cal X$, equipped
with norms $\|.\|_{t}$, such that
for all $t$ in $\R$, $X_{t}$ satisfies
 a $f$-lower estimate and
for all $x$ in $c_{00}$, the function $t \mapsto \|x\|_{t}$ is
 measurable.
For vectors of range at
most $E_N$, we have
 $f(N)^{-1}\|x\|_1 \leq \|x\|_t \leq \|x\|_1$, so
that the norms $\|.\|_t$ are equivalent to $\|.\|_1$
with log-integrable constants. We are then allowed
to define the $\theta$-interpolation space of the family defined
 on $\delta S$
as $X_{t}$ if $z=it$, $l_{q}$ if $z=1+it$.
 Let
${\cal X}_{\theta}$ be the class of  spaces $X$ obtained in
that way.

We shall sometimes use for $z \in \delta S$ the notation
$\|.\|_{z}$, to mean $\|.\|_{t}$ if $z=it$, and $\|.\|_{q}$ if
$z=1+it$. There will be no ambiguity from the context.
We shall similarly use the notation $\|.\|^{*}_{z}$. 
The notation $X_t^N$ stands for $E_N X_t$, and
$X_t^{N*}$ for $E_N X_t^*$. Also,
if not specified, the measure of a subset of $\R$ 
will be its measure for $\mu_0$.

\subsection{Properties of ${\cal X}_{\theta}$}
Let $X$ be in 
${\cal X}_{\theta}$ and $x$ be in $X$.
 Let ${\cal A}_{\theta} (x)$ be the set of
 functions in ${\cal A}_{S}$
 that take the value $x$ at the point $\theta$.
 Given $\theta$, it is the set of {\em interpolation 
functions for $x$}.
By definition, for all $x$ in $X$,
 $ \|x\| = \inf_{F \in {\cal A}_{\theta} (x)}
 (\int_{z \in \delta S} \|F(z)\|_{z}  d \mu (z))$.
The following theorem is a useful result of \cite {CC}.

\paragraph{Theorem 1}
{\em If $x$ is of finite range, there is an interpolation
function $F$ for $x$, that we shall call
 {\em minimal for $x$}, with
$ran(F)=ran(x)$ and such that} 
 
\[ \|F(it)\|_{t}=\|x\| \ a.e. \ and\ 
\|F(1+it)\|_{q}=\|x\| \ a.e.\ .\]

 \paragraph{Lemma 1}
{\em The following formula is also true:
 
 \[ \|x\| = \inf_{F \in {\cal A}_{\theta} (x)}
 \left(\int_{\R} \|F(it)\|_{t}  d \mu_{0}(t)\right)
^{1-\theta}
  \left(\int_{\R} \|F(1+it)\| _{q} d \mu_{1} (t)\right)
^{\theta} .\] }

{\em Proof } First notice that for any $F$ in ${\cal A}_{\theta} (x)$,
by a convexity inequality, the argument in the second infimum is 
smaller than
\[ (1-\theta) \left(\int_{\R} \|F(it)\|_{t}
  d \mu_{0}(t)\right) +
\theta \left(\int_{\R} \|F(1+it)\| _{q}
 d \mu_{1} (t)\right) \]
equal to
$\int_{z \in \delta S} \|F(z)\|_{z}  d \mu (z)$,
so that the second infimum is smaller than the first one.

Now, given $u \in \R$, the map $G_{u}$ defined on 
${\cal A}_{\theta} (x)$ by $ G_{u} (F) (z) = F(z) e^{u(z-\theta)} $
is a bijection on ${\cal A}_{\theta} (x)$. Furthermore, for any $u$,
the expressions 
\[ \left(\int_{\R} \|(G_{u}(F)(it)\|_{t} 
 d \mu_{0}(t)\right)^{1-\theta}
  \left(\int_{\R} \|G_{u}(F)(1+it)\| _{q} 
d \mu_{1} (t)\right)^{\theta} \]  and 
\[ \left(\int_{\R} \|F(it)\|_{t}  d \mu_{0}(t)\right)^{1-\theta}
  \left(\int_{\R} \|F(1+it)\| _{q} d \mu_{1} (t)\right)^{\theta} \]
are equal.
If we choose a proper $u$ (namely such that 
$\int_{\R} \|(G_{u}(F)(it)\|_{t}  d \mu_{0}(t) =
\int_{\R} \|G_{u}(F)(1+it)\| _{q} d \mu_{1} (t)$), this
is also equal to $\int_{z \in \delta S} \|G_{u}(F)(z)\|_{z}  d \mu (z)$.
Consequently, the two infima are actually equal.

\paragraph{Proposition 1}
{\em For all successive vectors
 $x_{1}<\cdots<x_{n}$ in $X$,
\[ \frac{1}{f(n)^{1-\theta}} 
\left( \sum_{i=1}^{n} \| x_{i} \| ^{p} \right)^{\frac{1}{p}}
\leq \Bigl\| \sum_{i=1}^{n} x_{i} \Bigl\|
\leq  \left( \sum_{i=1}^{n} \| x_{i} \| ^{p} \right)^{\frac{1}{p}} .\] }

{\em Proof } It is enough to prove this in the interpolation space
$X^{N}$ defined above, written in short
 $(X_{t}^{N},l_{q}^{N})_{\theta}$, for any $N \geq 1$.

\

{\em First inequality}
The unit ball of $X_{t}^{N}$ is stable under sums of the form 
$\sum_{j=1}^{n} \lambda_{j} y_{j}$, where the $y_{j}$ are successive
in the unit ball of $X_{t}^{N}$
and $\sum_{j=1}^{n}|\lambda_{j}| = 1$.

The unit ball of $l_{q}^{N}$ is stable under sums of the form 
$\sum_{j=1}^{n} \mu_{j} z_{j}$, where the $z_{j}$ are successive
in the unit ball of $l_{q}^{N}$and $\sum_{j=1}^{n}|\mu_{j}|^{q} = 1$.

Consequently, the unit ball of $X^{N}$ is sta\-ble under
successive sums
of the form\ 
 $\sum_{j=1}^{n} \lambda_{j}^{1-\theta} 
\mu_{j}^{\theta} x_{j} $, where the $x_j$ are in the unit ball of $X^N$ and
$\lambda_{j}$ and $\mu_{j}$ satisfy the above conditions. Indeed,
for every $x_{j}$ in the unit ball of $X^{N}$, let $F_{j}$ be
 minimal for $x_j$;
 the function $F$ defined by $F(z)=\sum_{j=1}^{n} \lambda_{j}^{1-z}
 \mu_{j}^{z} F_{j}(z)$ is then in ${\cal A}_{S}$ and bounded by $1\ a.e.$ on 
$\delta S$, so by definition, $\|F(\theta)\| \leq 1$, that is,
 $\sum_{j=1}^{n} \lambda_{j}^{1-\theta} \mu_{j}^{\theta} x_{j} $
is in the unit ball of $X^{N}$.

Now consider any successive vectors $x_{j}$ in $X^{N}$, and apply this
stability property to $x_{j}/ \|x_{j}\|$ and 
$\lambda_{j} = \mu_{j}^{q} =
 \|x_{j}\|^{p} / \sum_{i=1}^{n} \|x_{i}\|^{p} $. Using the equality
$1-\theta +\theta /q = 1/p$, one finally gets:
\[ \Bigl\| \sum_{j=1}^{n} x_{i} \Bigl\|
\leq  \left( \sum_{j=1}^{n} \| x_{i} \| ^{p} \right)^{\frac{1}{p}} .\]

This inequality will be called the {\em upper $p$-estimate for $X$}.

\

{\em Second inequality}
According to \cite {CC}, the duality property is true in finite
dimension, that is $X^{N*} = ( (X_{t}^{N})^{*},l_{q}^{N*} )_{\theta} $.
As $X_{t}$ satisfies a lower $f$-estimate, so does $X_{t}^{N}$; the
 dual version of this is that the unit ball of $(X_{t}^{N})^{*}$ is stable
under sums of the form
 $(1/f(n)) \sum_{j=1}^{n} y_{j}^{*}$,
 where the $y_{j}^{*}$ are successive. As $l_{q}^{N*}=l_{\qiq}^{N}$, we know
that its unit ball is stable under successive sums of the form 
$\sum_{j=1}^{n} \mu_{j} z_{j}^{*}$,
 where $\sum_{j=1}^{n}|\mu_{j}|^{\qiq} = 1$. Letting $\lambda_{j} = 
1 / f(n) $ for each $j$, and using the same
 proof as above, we get that
the unit ball of $X^{N*}$ is stable under
 successive sums of the form
$(1/f(n)^{1-\theta}) \sum_{j=1}^{n} \mu_{j}^{\theta} x_{j}^{*}$.

Now let $x_{j}$ be successive vectors in $X^{N}$; for $j=1, \ldots, n$,
let $x_{j}^{*}$ be successive dual unit vec\-tors such that
$x_{j}^{*}$ norms $x_{j}$ ( recall that the basis is
bimonotone in every $X_t$, so it is bimonotone in $X$). We get that
$ (1/f(n)^{1-\theta}) 
\sum_{j=1}^{n} \mu_{j}^{\theta} \|x_{j}\|
 \leq \| \sum_{j=1}^{n} x_{j} \| $.
Choosing $\mu_{j}^{\qiq} =
\|x_{j}\|^{p} / \sum_{i=1}^{n} \|x_{i}\|^{p}$ and using the
equality $\theta / \qq = 1-1/p$ gives the desired inequality:

\[ \frac{1}{f(n)^{1-\theta}} 
\left( \sum_{i=1}^{n} \| x_{i} \| ^{p} \right)^{\frac{1}{p}}
\leq \Bigl\| \sum_{i=1}^{n} x_{i} \Bigl\|. \]

This inequality will be called the {\em lower estimate for $X$}.

\

{\em Remark} Gowers-Maurey's space, and, more generally, spaces
satisfying $f$-lower estimates 'look like' the space $l_{1}$ ( for
successive vectors, the
triangular inequality is, up to a 
logarithmic term, an equality).
As the interpolation space of $l_{1}$ and $l_{q}$ is $l_{p}$, one
expects the space $X$ to 'look like' $l_{p}$; the above inequalities
show in what sense this is true.

\paragraph{Proposition 2}
{\em The dual space $X^{*}$ of $X$ is
 also the interpolation space - as defined at the end of
1.1 - of the family defined
 on $\delta S$
as $X_{t}^{*}$ if $z=it$ and $l_{\qiq}$ if $z=1+it$.}

{\em Proof } Recall that a basis $(x_{n})_{n=1}^{\infty}$ of
a Banach space is {\em shrinking} if for every continuous linear
functional $x^{*}$ and every $\epsilon > 0$ there exists
 $n \in \N$ such that the norm of $x^{*}$ restricted to the
span of $x_{n},x_{n+1},\ldots$ is at most $\epsilon$.
The basis $e_{1},e_{2},\ldots$ is a shrinking basis for $X$. Indeed,
 suppose it is not; then we can find
 $\epsilon > 0$, a norm-1 functional $x^{*} \in X^{*}$,
and a sequence of successive normalized blocks $x_{1},x_{2},\ldots$
such that $x^{*}(x_{n}) \geq \epsilon $ for every $n$. Then, using 
the upper $p$-estimate,
 we get $n \epsilon \leq x^{*}(\sum_{i=1}^{n} x_{i})
 \leq \| \sum_{i=1}^{n} x_{i} \| \leq n^{1/p} $,
a contradiction if we choose $n$ big enough.

This implies that given $x^{*}$ in $X^{*}$,
$\|x^{*}\|_{X^{*}} =
 \lim_{N \rightarrow +\infty} \| E_{N} x^{*} \|_{X^{N*}} $.
 But this means that 
$ X^{*} = completion (\cup_{n \in \N} X^{N*} ) $;
furthermore, according to \cite {CC}, $X^{N*}$ is
also the interpolation space
 $((X^N_t)^*,l_{\qiq}^N)_{\theta}$; as
$(X^N_t)^*=(X^*_t)^N$, we get
the desired dual property.

\paragraph{Proposition 3}
{\em The space $X$ is uniformly convex.}

{\em Proof } It is enough to prove that any  vectors $x$ and
$y$ in the unit ball of $X^N$ satisfy the relation
 $\| \frac{x+y}{2} \| \leq 1-\delta (\|x-y\|)$ where $\delta$ is
strictly positive on $]0,+\infty[$ and does not depend on $N$. 

We know by \cite {CC} that for any $r \geq 1$ the norm of a vector
$x$ in $X_{N}$ is given by the formula 
 $ \|x\|^{r} = \inf_{F \in {\cal A}_{\theta} (x)}
 (\int_{z \in \delta S} \|F(z)\|_{z}^{r}  d \mu (z))$.
 As in Lemma 1, we have also:

 \[ \|x\|^{r} = \inf_{F \in {\cal A}_{\theta} (x)}
 \left(\int_{\R} \|F(it)\|_{t}^{r}  
d \mu_{0}(t)\right)^{1-\theta}
 \left(\int_{\R} \|F(1+it)\|_{q}^{r}
 d \mu_{1}(t)\right)^{\theta} .\]

Suppose $q \geq 2$. Then for any vectors $a$ and $b$ in the
unit ball of $l_{q}^{N}$,
 $ \| \frac{a+b}{2} \|^{q}_{q} \leq 1 - \| \frac{a-b}{2} \|^{q}_{q} $
 ( this Clarkson's inequality can be found in \cite {B}).
 Now let $x$ and $y$ be in the unit ball of $X^{N}$, let
$F$ ( resp. $G$) be a minimal interpolation function for $x$ ( resp.
 $y$) as in Theorem 1. Let us apply the formula with $r=q$:
\[\left\| \frac{x+y}{2} \!\right\|^{q} \leq 
 \left(\int_{\R}
 \left\| \frac{F+G}{2}(it)\right\|_{t}^{q}  
d \mu_{0}(t)\right)^{1-\theta}
\!  \left(\int_{\R}
 \left\| \frac{F+G}{2} (1+it)\right\|_{q}^{q}
 d \mu_{1} (t)\right)^{\theta}\!\!.\]

The first integral is smaller than $1$, so that:
\[\left\| \frac{x+y}{2} \right\|^{q} \leq 
  \left(\int_{\R}
 \left\| \frac{F+G}{2} (1+it)\right\|_{q}^{q}
 d \mu_{1} (t)\right)^{\theta}. \]
Similarly,
\[\left\| \frac{x-y}{2} \right\|^{q} \leq 
  \left(\int_{\R}
 \left\| \frac{F-G}{2} (1+it)\right\|_{q}^{q}
 d \mu_{1} (t)\right)^{\theta}. \]
Adding  these two estimates together, and using Clarkson's estimate we get 
\[ \left\| \frac{x+y}{2} \right\| ^{q/\theta} +
\left\| \frac{x-y}{2} \right\| ^{q/\theta} \leq 1 .\]

If $q <2$, there is another estimate in \cite {B}:
 there is a constant $c_{q}$ such that
 for any vectors $a$ and $b$ in the
unit ball of $l_{q}^{N}$,
$ \| \frac{a+b}{2} \|_{q} \leq 1 - c_{q} \| a-b \|^{2}_{q} $.
 Applying the same method as above,
 we obtain

\[ \left\| \frac{x+y}{2} \right\| ^{1/\theta} + c_{q} \| x-y \| ^{2/\theta} \leq 1.\]

In both the cases $q \geq 2$ and $q<2$, the inequalities above
are uniform convexity inequalities.

\subsection{$l_{p+}^{n}$-averages}

\paragraph{Definition 1}

Let $n$ be a non-zero integer, $C$ a real number.

Let $X$ be in ${\cal X}$.
An {\em $l_{1+}^{n}$-average in $X$ with constant $C$}
 is a normalized vector $x \in X$ such that
$ x= \sum_{i=1}^{n} x_{i} $ where $x_{1}< \cdots < x_{n}$ are successive
 vectors and each $x_{i}$
 verifies $\|x_{i}\| \leq C n^{-1} $.

Let $X$ be in ${\cal X}_{\theta}$.
An {\em $l_{p+}^{n}$-average in $X$ with constant $C$}
 is a normalized vector $x \in X$ such that
$ x= \sum_{i=1}^{n} x_{i} $ where $x_{1}< \cdots < x_{n}$ are successive
 vectors and each $x_{i}$
 verifies $\|x_{i}\| \leq C n^{-1/p} $.

An {\em $l_{1+}^{n}$ ( resp. $l_{p+}^{n}$) -vector} is a non-zero
multiple of an $l_{1+}^{n}$ ( resp. $l_{p+}^{n}$)-average.

\paragraph{Lemma 2}
{\em Let $X$ be in ${\cal X}_{\theta}$. For every  $n \geq 1$,
 every $C >1$, every block subspace $Y$ of $X$ 
contains an $l_{p+}^{n}$-average with constant $C$.}

{\em Proof } The proof is the same as
in Lemma 3 of \cite{GM}. Suppose the result is false for some $Y$. Let $k$ be an integer
such that
 $ k \log C > (1-\theta) \log f(n^{k})  $,
 let $N=n^{k}$, let $x_{1} < \cdots < x_{N}$
be any sequence of successive norm-$1$ vectors in $Y$,
and let $x=\sum_{i=1}^{N} x_{i}$. For every $0 \leq i \leq k$
and every $1 \leq j \leq n^{k-i} $, let 
$x(i,j)= \sum_{t=(j-1)n^{i}+1}^{jn^{i}} x_{t}$. Thus
$x(0,j)=x_{j}, x(k,1)=x$, and, for $1 \leq i \leq k$, each
$x(i,j)$ is a sum of $n$ successive $x(i-1,j)$'s.
By our assumption, no $x(i,j)$ is an $l^{n}_{p+}$-vector
with constant $C$. It follows easily by induction that
$\| x(i,j) \| \leq C^{-i}n^{i/p} $ and, in particular, that
$\| x \| \leq C^{-k} n^{k/p} = C^{-k} N^{1/p} $. However, it 
follows from the lower estimate in $X$ that
$\| x \| \geq N^{1/p} f(N)^{-(1-\theta)} $. This is a 
contradiction, by choice of $k$.

\paragraph{Lemma 3}
{\em Let $X$ be in ${\cal X}_{\theta}$. Let $0<\epsilon<1/4$. Let 
$\theta = 1/2 $.
 Let $x$ be an $l_{p+}^{n}$-average in $X$  with constant
 $1+\epsilon$.
There exists an interpolation function $F$ for $x$ with
$ran(F)=ran(x)$, bounded almost everywhere by $1+\epsilon$,
such that except on a set of measure at most
$2 \sqrt{\epsilon}$, $F(it)$ is an $l_{1+}^{n}$-vector in $X_{t}$, of norm $1$
 up to $\sqrt{\epsilon}$,
 with constant $1+ 4 \sqrt{\epsilon}$.}

 Such a function is called
 {\em $\epsilon$-representative}, or {\em representative}, since
we shall always consider $l^n_{p+}$-averages associated
to given values of $\epsilon$.

{\em Proof} The vector $x$ can be written $\sum_{j=1}^{n}x_{j}$
where $x_{1}< \cdots < x_{n}$ are successive
 vectors and each $x_{j}$
 verifies $\|x_{j}\| \leq (1+\epsilon) n^{-1/p} $.
Let $F^{\prime}_{j}$ be a minimal interpolation function for $x_{j}$,
let $F_{j}$ be defined by $F_{j}(z) = n^{-1/ \pip + z/ \qiq}
 F^{\prime}_{j} (z)$
and let $F=
 \sum_{j=1}^{n} F_{j} $. We show
that $F$ is representative for $x$.

 Notice that $F(\theta)=x$, so 
$F$ is an interpolation function for $x$, and

\[ 1= \|x\| \leq 
 \left(\int_{\R} \|F(it)\|_{t}   d \mu_{0}(t)\right)^{1-\theta}
  \left(\int_{\R} \|F(1+it)\| _{q} 
d \mu_{1} (t)\right)^{\theta}.\]

By choice of $F$, $F$ is bounded by
$1+\epsilon\ a.e.$ on $\delta S$, so both integrals are smaller than $1+\epsilon$.
 As a consequence, 
$\int_{t \in \R} \|F(it)\|_{t}  d \mu_{0}(t)
\geq (1+\epsilon)^{-\theta / (1-\theta)} \geq 1- \epsilon $ ( recall
that $\theta = 1/2$). As for every $t$,
 $\| F(it) \|_{t} \leq 1+\epsilon$, by a Bienaym\'{e}-Tchebitschev
estimation, we get that on a set of measure at least 
$1-2 \sqrt{\epsilon}$, $\| F(it) \|_{t} \geq 1-\sqrt{\epsilon}$.

So on that set, $F(it)$ is of norm $1$ up to $\sqrt{\epsilon}$.
 For each $j$,
$ \|F_{j}(it)\|_{t} = n^{-1/ \pip} \|x_{j}\| \leq (1+\epsilon) / n $;
 so that $F(it)$ is an 
$l_{1+}^{n}$-vector in $X_{t}$ with constant 
$(1+\epsilon)/(1-\sqrt{\epsilon}) \leq 1+4 \sqrt{\epsilon}$.

\subsection{ Rapidly Increasing Sequences }

To make a construction similar to the one in \cite{GM}, we need definitions
of Rapidly
Increasing Sequences in $X$ and of special sequences in $X^{*}$. 
We now assume that $\theta = 1/2$.

\paragraph{Definition 2}

Let $N$ be a non-zero integer. Let $0<\epsilon \leq 1$.

Let $X$ be in ${\cal X}_{\theta}$.
A sequence $x_{1} < \cdots < x_{N}$ in $X$ is a
 {\em Rapidly Increasing Sequence of $l_{p+}^{n}$-averages,
 or R.I.S., of length $N$ with constant $1+\epsilon$} if $x_{k}$ is an 
$ l_{p+}^{n_{k}}$-average with constant $1+\epsilon / n_{k}$ for each $k$,
$ n_{1} \geq 4 M_{f} (N/ \epsilon )/ \epsilon \ff (1) $,
 and $\epsilon /2 \ f(n_{k})^{1/2} \geq |ran(x_{k-1})| $ 
for $k=2,\dots,N$.

Here $f^{\prime}(1)$ is the right derivative of $f$ at $1$
and $M_f$ is defined on $[1,\infty)$ by $M_f(x)=f^{-1}(36x^2)$.
 
 In spaces $X_{t}$, we shall use R.I.S. in Gowers-Maurey sense,
 that is,
sequences of 
$ l_{1+}^{n_{k}}$-averages with constant $1+\epsilon$ with the
same increasing condition as above.

We shall call both kinds "R.I.S." without ambiguity.
A R.I.S.-vector is a non-zero multiple of the sum of a R.I.S..
The following proposition links the two kinds of R.I.S..

\paragraph{Lemma 4}
{\em Let $X$ be in ${\cal X}_{\theta}$. Let $0 < \epsilon \leq 1/16$.
Let $x_{1} < \cdots < x_{n}$ be a R.I.S. in $X$
 with constant $1+\epsilon$, and let $x=\sum_{k=1}^{n}x_{k}$.
For each $k$, let $F_{k}$ be representative for $x_{k}$;
 then 
$F = F_{1}+ \cdots + F_{n}$ 
is an interpolation function for $x$, and
except on a set of measure at most
$4 \sqrt{\epsilon} / f(n)$, $F(it)$ is up to $ 2 \sqrt{\epsilon}$
 the sum of a
R.I.S. in $X_{t}$ with constant $1+ 4 \sqrt{\epsilon}$.}

{\em Proof } It is clear that $F$ is an interpolation
 function for $x$. According to Lemma 3, for 
each $k$, $F_{k}(it)$ is 'close' to an $l_{1+}^{n_{k}}$-average,
except on a set of measure at most $2 \sqrt{\epsilon/n_{k}}$.
The union over $k$ of these sets is of measure at most
$\sum_{k=1}^{n} 2 \sqrt{\epsilon / n_{k}} 
\leq 4 \sqrt{\epsilon / n_{1}}
\leq 4 \sqrt{\epsilon} / f(n) $
 ( this is a consequence of
the increasing condition and the lower bound for $n_1$ in the definition of the R.I.S.).

Now let $t$ be in this union.
 For every $k$, let $|F|_{k}(it)$ denote
the normalization of $F_{k}(it)$; $|F|_{k}(it)$ is
 an $l_{1+}^{n_{k}}$-average
 with constant $1+4 \sqrt{\epsilon / n_{k}}$.
The sequence $|F|_{1}(it) < \cdots < |F|_{n}(it) $ is
 a R.I.S. in $X_{t}$, with constant
 $\sup_{k} (1+4 \sqrt{\epsilon / n_{k}}) \leq 1+4\sqrt{\epsilon}$
( because $1+4\sqrt{\epsilon} > 1+\epsilon$, the increasing condition is indeed
verified).

It remains to show that $F(it)$ and the sum of the $|F|_{k}(it)$ are
equal up to $2 \sqrt{\epsilon} $; and indeed
$ \| F(it)- \sum_{k=1}^{n} |F|_{k}(it) \|_{t}
\leq \sum_{k=1}^{n} | 1-\|F_{k}(it)\| |_{t}
\leq \sum_{k=1}^{n} \sqrt{\epsilon / n_{k}}
 \leq 2 \sqrt{\epsilon}$,
so that the proof is complete.

\paragraph{Special sequences}

The trick is to define special sequences of {\em dual interpolation functions}.
 Thus, by a Gowers-Maurey construction, we obtain
spaces $X_{t}$ that "look like" Gowers-Maurey's
space and such that
the special property of the $X_{t}$ is somehow
uniform on $t$; more precisely, we build 
spaces $X_{t}$ - and the related $X$ - and a space $\Delta$ of dual
interpolation functions
 such that
 $\Delta$ is countable,
stable under 'Schlumprecht's operation' and under taking special 
functions, and such that 
any vector in the unit ball of $X^{*}$ has an almost
minimal interpolation function in $\Delta$.
This construction, and the proof that $X$ is
hereditarily indecomposable, are developped in the next two parts.

\section{ Construction of a space $X$ in ${\cal X}_{\theta}$ }

\subsection{ Construction of spaces $X_{t}$}

Let $J=\{j_{1},j_{2},\ldots\}$, where $(j_{n})_{n \in \N}$ is a sequence
of integers such that $f(j_{1})>256$ and $\log \log \log j_{n}>4(j_{n-1})^{2}$
 for $n>1$.
Let $K=\{j_{1},j_{3},j_{5},\ldots \}$ and $L=\{j_{2},j_{4},j_{6},\ldots \}$. 
Let  $\{ L_{m}, m \in \N^* \}$ be a partition of $L$ with every $L_{m}$ infinite.
For $r \in [1,+\infty] $, let $B(l_{r})$ denote the
 unit ball of  $l_{r} \cap c_{00}$.
For $N \geq 1$ and $z \in \C$,
 let $f(N,z)=f(N)^{1-z}N^{z/\qiq}$ and
$g(N,z)=\sqrt{f(N)}^{1-z}N^{z/\qiq}$.

\paragraph{Definition 3}
Given a subset $D$ of ${\cal A}_{S}$,
for every $N>0$, the set of
 {\em $N$-Schlum\-precht sums in $D$}, written $B_{N}(D)$, is the set
of functions of the form $f(N,z)^{-1} 
 \Sigma_{i=1}^{N} F_{i}$, where the $F_i$ are successive in $D$.
A {\em Schlumprecht sum in $D$} is a $N$-Schlumprecht sum in $D$
for some
$N > 0$.
Let $B(D)$ be the set of Schlumprecht sums.
If $D$ is countable, given an injection $\tau$
 from $ \cup_{m \in \N}B(D)^{m} $
 to $\N$, and an integer $k$, a {\em special function} in $D$, for $\tau$, with length $k$, is a
function of the form
  $g(k,z)^{-1}  \Sigma_{j=1}^{k} G_{j}$, with 
 $ G_{j} \in B_{n_{j}} (D), G_{1}< \cdots <G_{k}$, $n_{1}=j_{2k}$ and $ n_{j}=
 \tau (G_{1},\ldots,G_{j-1}) $ for $j=2,\ldots,k$;
 $G_{1},\ldots,G_{k}$ is 
a {\em special sequence in $D$}.

Here, it does not seem possible to define the set of special functions 
before defining the spaces $X_t$
as in \cite{GM}, so we build them at the same time by induction.

\paragraph{Step 1}

For every $t$ in $\R$, let $D_{1}(t) = B(l_1)$. 
Let ${\cal D}_{1}$ be the set of functions
 in ${\cal A}_{S}$  with values
 in $D_{1}(t)$ for almost every $it$ and in $B(l_{\qiq})$ almost
 everywhere
 on $S_{1}$.
Let $\Delta_{1}$ be a countable set of  
functions in ${\cal A}_{S}$, dense in
 ${\cal D}_{1}$ for the $L_{1}$-norm ( namely
 $ \|F\|= \int_{z \in \delta S} \|F(z)\|_{1}  d \mu (z) $). For
this first step, we may assume that all functions in $\Delta_1$
are continuous.
Let $\sigma_{1}$ be an injection from
 $ \cup_{m \in \N} (\Delta_{1})^{m} $ to 
$L_{1} $, the first subset of $L$ in the partition
 $\{ L_{m}, m \in \N^* \}$. Let $S_0^1=\R$.

\paragraph{Step n}
We are given a set of sequences $D_{n-1}(t)$ for every $t$ in $\R$,
 a set ${\cal D}_{n-1}$ of functions in ${\cal A}_S$, a countable set
 $\Delta_{n-1}$ of functions in ${\cal A}_S$ defined
 everywhere on $S_0$,
a subset $S_0^{n-1}$ of $\R$ of measure $1$ ( that
stands for the set of 'significative' values of the functions
in $\Delta_{n-1}$), and an
 injection $\sigma_{n-1}$ from 
$ \cup_{m \in \N} (\Delta_{n-1})^{m} $ to 
$L_{1} \cup \ldots \cup L_{2n-3} $.

Then let $ \DeDe_{n} =B(\Delta_{n-1})
  \cup \{ EF, E \ interval, F \in \Delta_{n-1} \}
 $.
Let $\tau_{n-1}$ be an injection from $ \cup_{m \in \N} 
(\DeDe_{n}^{m} \setminus \Delta_{n-1}^{m}) $ to 
$L_{2n-2} $.

Let ${\cal S}_{n-1}$ be the set of special functions in $\Delta_{n-1}$,
 for $\tau_{n-1}$, with length in $K$.
For every $t$ in $\R$,
let $D_{n}^{S}(t)$
be the sets of sequences of the form
$f(N)^{-1} \Sigma_{i=1}^{N} x_{i}$ where
the $x_i$ are successive in $D_{n-1}(t)$,
$D_{n}^{I}(t)$ be the set of sequences $Ex$ where $E$ is an interval
and $x$ is in $D_{n-1}(t)$;
if $t \in S_0^{n-1}$, let 
$D_n^{s}(t)$ be the set of sequences of the form $G(it)$
where $G \in {\cal S}_{n-1}  $, otherwise, let $D_n^{s}(t)=\emptyset$.
 Let 
$D^{\prime}_n(t)=D_{n}^{S}(t) \cup D_{n}^{I}(t) \cup D_{n}^{s}(t)$
and
let $D_{n}(t)= conv ( \cup_{|\lambda|=1} \lambda \DD_{n}(t))$.
Let ${\cal D}_{n}$ be the set of functions
 in ${\cal A}_{S}$ with values in $D_{n}(t)$
 for almost every $it$ and in $B(l_{\qiq})$ almost everywhere on $S_{1}$.

We complete $\DeDe_{n}$ in $\Delta_{n}$ countable set of functions
in ${\cal A}_{S}$, dense in ${\cal D}_{n}$ for the $L_1$-norm. 
There is a subset $S_0^n \subset S_0^{n-1}$
 of $\R$ of measure $1$ such that
$F(it)$ is indeed in $D_n(t)$ for all $F$ in $\Delta_n$ and
for all $t$ in $S_0^n$.
 With an 
injection ${\tau}^{\prime}_{n-1}$, from 
$ \cup_{m \in \N} (\Delta_{n}^{m} \setminus \DeDe_{n}^{m}) $ to 
$L_{2n-1} $, we obtain an injection $\sigma_{n}$,
 from $ \cup_{m \in \N} (\Delta_{n})^{m} $ to 
$L_{1} \cup \ldots \cup L_{2n-1} $.

\paragraph{Definition of $X_{t}$}
It is easy to verify that the sequences $D_{n}(t)$ for every $t$
 in $\R$, ${\cal D}_{n}$ and
 $\Delta_{n}$ are increasing, that
the sequence $S_0^n$ is decreasing and that for every $n$, $\sigma_{n}$ coincides with $\sigma_{n-1}$
 on its set of definition.

 We then define  $D_{t}=\cup_{n \in \N} D_{n}(t)$ for every $t$ in $\R$,
${\cal D} =\cup_{n \in \N}{\cal D}_{n}$, $\Delta=\cup_{n \in \N}
\Delta_{n}$, $S_0^{\infty}=\cap_{n \in \N}S_0^{n}$ and
 $\sigma$ the injection from $ \cup_{m \in \N} \Delta^{m} $
 to 
$L$ whose restrictions are the $\sigma_{n}$.

Finally for every $t$ in $\R$, we define the space
 $X_{t}$ by its norm on $c_{00}$:
\[ \forall x \in c_{00}, \|x\|_t= \sup_{y \in D(t)} | <x,y>| .\]

\subsection{Properties of $\cal D$ and $\Delta$}

\paragraph{Proposition 4}\ 

{\em
(a) For every $t$ in $\R$,
 $B(l_{1}) \subset D(t) \subset B(l_{\infty})$.

(b) The set $\Delta$ is  countable, dense in $\cal D$, stable under
 interval projections and Schlumprecht sums in $\Delta$.

(c) For every $t$ in $\R$, the set $D(t)$ is convex,
 stable under interval
projections, multiplication by a scalar of modulus $1$ 
( or {\em balanced}), and
 sums of the form $ f(N)^{-1}  
\Sigma_{i=1}^{N} x_{i}$, with $ x_{i} \in D(t)$
 and $x_{1}< \cdots <x_{N} $.

(d) The set $\cal D$ is  convex, balanced, stable under
 interval projections, Schlum\-precht sums in $\cal D$, and under
 taking special functions in $\Delta$ for $\sigma$ with length in $K$.
}

{\em Proof}

(a) The left inclusion is a consequence of the facts that
$D_{1}(t) = B(l_{1})$ and that $D_{n}(t)$ is increasing;
for the right inclusion, notice that by induction,
 $D_{n}(t) \subset B(l_{\infty})$ for every $n$.

(b) The set $\Delta$ is countable as a countable union of countable
sets; it is dense in $\cal D$ because for every $n$,
$\Delta_{n}$ is dense in ${\cal D}_{n}$; the stability property under
interval projections and Schlumprecht sums is ensured because
for every $n$, $\Delta_{n}$ contains $\DeDe_{n}$, the set 
of projections and  sums from $\Delta_{n-1}$.

(c) The set $D(t)$ is convex as an increasing
 union of convex sets; the stability properties are ensured
by the definition of $\DD_{n}(t)$ and $D_{n}(t)$ from $D_{n-1}(t)$.

(d) The set ${\cal D}_{n}$ is the set of functions with values
in the convex, balanced, and interval projection
 stable sets $D_{n}(t)$ and $B(l_{\qiq})$ on $\delta S$, so
that it is convex, balanced and stable
 under interval projections; and so is $\cal D$.
 
To show the Schlumprecht stability property, it is enough,
 given successive functions $F_1<\cdots<F_N$ in ${\cal D}_{n-1}$, to show
that $F=f(N,z)^{-1}  
\Sigma_{j=1}^{N} F_{j}$ is in ${\cal D}_{n}$.
For each $j$, $F_{j}(it)$ is in $D_{n-1}(t)$ almost everywhere.
 The set of  
$t \in \R$ such that this happens for every $j$ is still of
 measure $1$.
On this set,
 $F(it)= (f(N) N^{-1/ \qiq})^{it} (1/f(N)) \sum_{j=1}^{N} F_{j}(it)$
 is in
$D_{n}(t)$, by the definition of $D_{n}(t)$. In the same way, almost everywhere on $S_{1}$, $F_{j}(1+it)$
is in $B(l_{\qiq})$ for every $j$, so that
 $F(1+it)$ is in $B(l_{\qiq})$ too. By definition, this means that
$F$ is in ${\cal D}_{n}$.

To show the special property, first notice that a special function $G$
in $\Delta$ is a special function in $\Delta_{n}$ for some $n$
in $\N$.
It follows that $G(it) \in D_{n+1}(t)$ for every $t$ in $S_0^n$, that
is, almost everywhere;
furthermore,
 $G(1+it)$ is in $B(l_{\qiq})$ almost everywhere; so $G$
is in ${\cal D}_{n+1}$.

\paragraph{Lemma 5}
{\em Let $\cal S$ be the set of functions in ${\cal A}_{S}$ with values
 in $D(t)$ for almost every $it$ and in $B(l_{\qiq})$ almost
 everywhere on $S_{1}$. Then $\cal D$
 is dense in ${\cal S}$ for the $L_{1}$-norm.}

{\em Proof } Let $F$ be in ${\cal S}$, $ 0 <\epsilon < 1$.
 Let $N$ be such
that $ran(F) \subset E_{N}$.

 We recall Havin lemma from \cite{P} in a rougher version. Furthermore,
we state it on $S$ instead of on the unit disk of $\C$ ( the
two versions are equivalent using a conformal mapping ).

{\em Lemma } For every $\epsilon^{\prime} > 0$, there exists $\delta > 0$ such 
that for every subset $e$ of $\delta S$ with $\mu(e) \leq \delta $,
there exists $g_e$ in $H^{\infty}(S)$ with 
$|g_e| \leq 1 \ a.e. $ on $\delta S$,
 $\sup_{z \in e}|g_e (z)| \leq \epsilon^{\prime}$, and
$\int_{\delta S} |g_e (z) -1|d\mu(z) \leq \epsilon^{\prime}$.

Now let $\delta$ be
associated to $\epsilon^{\prime}=\epsilon / N$.
The sequence $( \{ t : F(it) \in D_{n}(t) \} )_{n \in \N}$ is
increasing and its union is of measure 1 for
 $\mu_{0}$, so there exists $n$ such that 
$T=\{ t / F(it) \in D_{n}(t) \}$ is of measure at least $1-\delta$.
For $\mu$, $\delta S \setminus iT$ is of measure at most
$\delta(1-\theta) \leq \delta$.
Let $H$ be the function $g_{\delta S \setminus iT}$.
 Let $\tilde{F}=H.F$. The function
$\tilde{F}$
is in ${\cal A}_{S}$. Furthermore,
 $\tilde{F}(1+it)$ is in $B(l_{\qiq})\ a.e.$
 on $S_{1}$, $\tilde{F}(it)$ is in $D_{n}(t)\ a.e.$ on $T$; this last
assertion
is also true on $S_{0} \setminus T$, because almost
everywhere on this set,
$ \tilde{F}(it)$ is in $1/N \ D(t)$ and because, for
functions of range at most $E_{N}$, we have
the following inclusions: $ 1/N \ D(t) \subset 1/N \ B(l_{\infty})
\subset B(l_{1}) \subset D_{n}(t)$. This proves that $\tilde{F}$
is in ${\cal D}_{n}$.

It remains to show that $F$ and $\tilde{F}$
are close, and indeed:
\[ \int_{\delta S}\|(F-\tilde{F})(z)\|_{1}d\mu(z)
\leq N \int_{\delta S}|H(z)-1|d\mu(z) \leq \epsilon. \]

\subsection{Definition of $X$}
\paragraph{}
 For every $x$ in $c_{00}$, the function $t \mapsto \| x \|_{t}$
is measurable. To see it,
it is enough to prove that the restriction of the
function to $S_0^{\infty}$ is measurable.
 We prove this by induction on $|ran(x)|$.
Remember that $\|x\|_{t}=\sup_{y \in D(t)}|<x,y>|$.
Now let $y$ be in $D(t)$;
 either $y$
is, up to multiplication by a scalar
of modulus $1$, the value in $it$ of the projection of
 a special function, and there are
countably many of them; or $y$ is a $n$-Schlumprecht sum with $n>1$
so that
$|<x,y>| \leq (1/f(n)) \sum_{j=1}^{n}\|{\cal E}_{j}x\|_{t} $, where
${\cal E}_1 < \cdots < {\cal E}_n$ are successive intervals; or $y$ is in $B(l_1)$.
Finally,

\[ \|x\|_{t}=\|x\|_{\infty}
 \bigvee \sup_{G \ special,E}|<x,EG(it)>|\ 
 \bigvee \sup_{n\geq2, {\cal E}_1 < \cdots < {\cal E}_n} \frac{1}{f(n)}
\sum_{j=1}^{n}\|{\cal E}_{j}x\|_{t}. \]

We may restrict the last sup to intervals ${\cal E}_{j}$ that do not
contain $ran(x)$; $t \mapsto \|x\|_{t}$ is then the supremum of
a countable family of measurable functions by the induction
hypothesis, so it is a measurable function.

 Furthermore,
it follows from the stability property of $D(t)$ that 
for every $t$ in $\R$, $X_{t}$ satisfies a lower-$f$
estimate.
We can then define a Banach space $X$ in ${\cal X}_{\theta}$ as in
 the first part of this article.

\paragraph{Lemma 6} {\em Let $F^{*} \in {\cal D}$. 
 Then $F^{*}(\theta)$ is in the unit ball of $X^{*}$.}

{\em Proof } First notice that if we
restrict them to finite range vectors, it is a consequence of
their convexity and of
the definition of $\|.\|_{t}$  that the unit
ball of $X_{t}^{*}$ and $\overline{D(t)}$ coincide.
Now, given $F^{*}$ in ${\cal D}$, it is of
 finite range. For almost every~$t$,
$F^{*}(it) \in D(t)$, so that by the previous remark,
$\| F^{*}(it) \|^{*}_{t} \leq 1$. Furthermore,
$\| F^{*}(1+it) \|_{\qiq} \leq 1$, so by Proposition 2,
$\| F^{*}(\theta) \|^* \leq 1$.

\

We 
need to recall some definitions and properties of \cite {GM}.
 Let $\cal F$ be Schlum-precht's space
of functions ( the explicit definition is in \cite {GM}; just
think of these functions as log-like).
We notice that $f$ and $\sqrt{f} \in {\cal F}$.
Given $X$ in $\cal X$,
 given $g$ in $\cal F$, a functional $x^{*}$ in $X^*$
is an $(M,g)-form$ if
$\|x^{*}\|^* \leq 1$ and $x^{*}= \sum_{j=1}^{M} x_{j}^{*}$ for
some sequence $x_{1}^{*} < \cdots < x_{M}^{*}$ of successive
functionals such that $\|x_{j}^{*}\|^* \leq g(M)^{-1}$ for each $j$.

Let $K_{0} \subset K$, and let us define
a function $\phi:[1,\infty) \mapsto [1,\infty)$ as
\[ \phi(x) = \sqrt{f(x)}\ \hbox{\ if}\ x \in K_{0},
   \ \phi(x) = f(x) \    \hbox{\ otherwise}. \]

We now state two lemmas that are a slight
modification of Lemma 7 of \cite {GM} for the first one and a mixture of 
Lemmas 8 and 9 of \cite {GM} for the second one.

Then we prove that the property of minimality
of the R.I.S. ( Lemma 10 of \cite {GM})
 is true in every $X_{t}$, and then that it can be extended to $X$.

\paragraph{Lemma 7}
{\em Let $f,g \in {\cal F}$ with
$g \geq \sqrt{f}$, let $X \in {\cal X}$ satisfy a lower $f$-estimate,
let $0<\epsilon \leq 1 $, let $x_{1} < \cdots < x_{N}$ be a R.I.S. in $X$
for $f$ with constant $1+\epsilon$, and 
let $x=\sum_{i=1}^{N} x_{i}$. Suppose that
\[ \|Ex\| \leq 1 \vee \sup \{|x^{*}(Ex)|: M \geq 2,
 x^{*}\ is\ an\ (M,g)-form\} \]
for every interval $E$. Then $\|x\| \leq (1+2\epsilon)
N g(N)^{-1}$.}

\paragraph{Lemma 8} {\em Given $K_{0} \subset K$, there is a function
 $g:[1,\infty) \mapsto [1,\infty)$
such that: $g \in {\cal F}$, $\sqrt{f} \leq g \leq \phi \leq f$,
and if $N \in J \setminus K_{0}$, then $g=f$ on the interval
 $[ \log N, \exp N ]$.}

\paragraph{Lemma 9}
{\em Let $t \in \R$.
Let $N \in L$, let $n \in [\log N,\exp N]$, let $0<\epsilon \leq 1$,
 and
let $x_{1}< \cdots <x_{n}$ be a R.I.S. in $X_{t}$ with constant
 $1+\epsilon$. Then 
\[ \Bigl\| \sum_{i=1}^{n} x_{i} \Bigl\|_{t} \leq (1+2 \epsilon) n / f(n). \] }

{\em Proof } 
The space $X_{t} \in {\cal X}$ satisfies a lower $f$-estimate.
 
 Let $x$ be the sum of the R.I.S.\ $x_{1}<\cdots<x_{n}$.
 Let $E$ be any interval.
 Let $\phi$ be the function
defined above in the case $K_{0}=K$ and $g$ associated to $\phi$
by Lemma 8. 
Let $x^{*}$ be a functional in $D(t)$. If $x^{*}$ is 
in $D_{1}(t)$, then
$|x^{*}(Ex)| \leq 1$. Else there exists $m \geq 2$ such that $x^{*}$ is
in $D_{m}(t) \setminus D_{m-1}(t)$; then, by definition of $D_{m}(t)$,
either $x^{*}$ is an $(M,f)-form$ with $M \geq 2$ 
or $x^{*}$ is an $(M,\sqrt{f})-form$
with $M \in K$; since $g \leq \phi$, it follows that $x^{*}$ is an
$(M,g)-form$ with $M \geq 2$. Consequently, 
\[ \|Ex\|_{t} \leq 1 \vee \sup\{|x^{*}(Ex)|:
 M \geq 2,\ x^{*}\ is\ an\ (M,g)-form\} \]
Since $g \in {\cal F}$ and $g \geq \sqrt{f}$, all the hypotheses
of Lemma 7 are satisfied. It follows that 
$\| \sum_{i=1}^{n} x_{i} \|_{t} \leq (1+2 \epsilon) n/g(n)$.
By Lemma 8, $g(n)=f(n)$, which proves our statement.

\paragraph{Lemma 10}
{\em Let $X$ be the space defined at
the beginning of 2.3. Let $N \in L$, let $n \in [\log N,\exp N]$,
 let $0 < \epsilon \leq 1/16$, and
let $X_{1}< \cdots <X_{n}$ be a R.I.S. in $X$ with 
constant $1+\epsilon$. Then
\[ \Bigl\| \sum_{i=1}^{n} X_{i} \Bigl\| \leq
 (1+10 \sqrt{\epsilon}) n^{1/p} / f(n)^{1-\theta}. \] }

{\em Proof }
Let
$F_{k}$ be representative for $X_{k}$, and $F=F_{1}+\ldots+F_{n}$.
We know that $F$ is an interpolation
 function for $X_1+\cdots+X_k$ so

 \[ \Bigl\|\sum_{i=1}^{n} X_{i}\Bigl\| \leq 
 \left(\int_{\R} \|F(it)\|_{t}  d \mu_{0}(t)\right)
^{1-\theta}
  \left(\int_{\R} \|F(1+it)\| _{q} d \mu_{1} (t)\right)
^{\theta} .\]

For the second integral, the following estimate holds:
\[ \int_{\R} \|F(1+it)\| _{q} 
d \mu_{1} (t) \leq (1+\epsilon) n^{1/q}. \]
According to Lemma 4, there is a set $A$ of
measure at most $4 \sqrt{\epsilon} / f(n) $ such that on
$\R \setminus A$, $F(it)$ is up to $2 \sqrt{\epsilon} $
 the sum $x_{t}$ of
a R.I.S. in $X_{t}$. On $\R \setminus A$,
$\| F(it) \|_{t} \leq \| x_{t} \|_{t} + 2 \sqrt{\epsilon} $;
furthermore, $x_{t}$ is a R.I.S. in $X_{t}$ with constant 
$1+4\sqrt{\epsilon}$, so that by Lemma 9,
$\| x_{t} \|_{t} \leq (1+8 \sqrt{\epsilon}) n/f(n) $.
On $A$, we have only $\| F(it) \|_{t} \leq (1+\epsilon) n.$

Gathering these estimates, we get:
\[ \int_{\R} \|F(it)\|_{t}  d \mu_{0}(t)
\leq [(1+8 \sqrt{\epsilon}) \frac{n}{f(n)} + 2 \sqrt{\epsilon} ] +
\frac{4 \sqrt{\epsilon}}{f(n)}\ (1+\epsilon) n \leq
(1+15 \sqrt{\epsilon}) \frac{n}{f(n)}.\]

Going back to the R.I.S. $X_{1} < \cdots <X_{n}$, we have
\[ \Bigl\| \sum_{i=1}^{n} X_{i} \Bigl\| \leq (1+15 \sqrt{\epsilon})^{1-\theta}
(1+\epsilon)^{\theta} \frac{n^{1-\theta+\theta / q}}{f(n)^{1-\theta}}
\leq (1+10 \sqrt{\epsilon}) \frac{n^{1/p}}{f(n)^{1-\theta}}. \]

\paragraph{Lemma 11}
{\em Let $t \in \R$.
Let $N\in L$, let $0<\epsilon<1/4$, let $M=N^{\epsilon}$
and let $x_{1}<\cdots<x_{N}$ be a R.I.S. in $X_{t}$
with constant $1+\epsilon$. Then $\sum_{i=1}^N x_i$ is an
$l_{1+}^M$-vector in $X_{t}$ with constant $1+4\epsilon$.}

{\em Proof }  It is the same as the one of Lemma 11
in \cite{GM}. Let $m=N/M$, let $x=\sum_{i=1}^N x_i$
 and for $1\leq j\leq M$ let
$y_j=\sum_{i=(j-1)m+1}^{jm}x_i$. Then each $y_j$ is
 the sum of a R.I.S. of
length $m$ with constant $(1+\epsilon)$. By Lemma 9 we have
$\|y_j\|_{t} \leq (1+2\epsilon) m f(m)^{-1}$ for every $j$ while
 $\| \sum_{j=1}^m y_j\|_{t}=\|x\|
\geq Nf(N)^{-1}$. It follows that $x$
is an $\l_{1+}^M$-vector in $X_{t}$ with constant at
most $(1+2\epsilon)f(N)/f(m)$. But $m=N^{1-\epsilon}$ so
 $f(N)/f(m) \leq (1-\epsilon)^{-1}$. The
result follows.

\paragraph{Lemma 12}
{\em Let $\epsilon_{0} = 1/10 $.
 Let $k \in K$ and $F_{1}^{*},\ldots,F_{k}^{*}$ be a special
sequence of length $k$, with $F_{i}^{*} \in B_{M_{i}}(\Delta)$.
 Let $t \in S_0^{\infty}$.
Let $x_{1}<\cdots<x_{k}$ a sequence of successive vectors in $X_{t}$,
 where every $x_{i}$ is a normalized
R.I.S.-vector of length $M_{i}$ and
 constant $1+\epsilon_{0}/4$. Suppose
 $ran(F_{i}^{*}) \subset
ran(x_{i})$ for $i=1,\ldots,k$,
 and $ 1/2 \ \epsilon_{0} f(M_{i}^{\epsilon_{0}/4})^{1/2}
 \geq |ran(x_{i-1})| $ for $i=2,\ldots,k$.

 If for every interval $E$, 
 $|(\sum_{i=1}^{k} F_{i}^{*}(t))(\sum_{i=1}^{k} Ex_{i})| \leq 4 $,
 then
\[ \Bigl\| \sum_{i=1}^{k} x_{i} \Bigl\|_{t} 
\leq (1+2 \epsilon_{0}) k / f(k) . \] }

{\em Proof }
First we recall two lemmas of \cite {GM}.

{\em Lemma GM4 } Let $M,N \in \N$ and $C\geq 1$, let $X\in \cal X$,
 let $x\in
X$ be an $l_{1+}^N$-vector with constant $C$ and
 let ${\cal E}_{1}<\cdots<{\cal E}_{M}$ be a
sequence of intervals. Then
$\sum_{j=1}^M \|{\cal E}_j x\| \leq C(1+2M/N) \|x\| .$

{\em Lemma GM5 } Let $f,g\in\cal F$ with $g\geq f^{1/2}$ and
 let $X\in\cal X$
satisfy a lower $f$-estimate. Let $0<\epsilon \leq 1$,
 let $x_1<\cdots<x_N$ be a R.I.S. in $X$
with constant $1+\epsilon$ and let $x=\sum_{i=1}^N x_i$.
 Let $M\geq M_f(N/\epsilon)$,
let $x^*$ be an $(M,g)$-form and let $E$ be any interval.
 Then $|x^*(Ex)|\leq
1+2\epsilon$.

\

 According to Lemma 11, each $x_{i}$ is an $l_{1+}^{N_i}$-average
with constant $1+\epsilon_{0}$, where $N_i = M_i^{\epsilon_0 /4}$.
 The increasing condition and the lower bound for
$M_1$ ensure that $x_1 < \cdots < x_k $ is a R.I.S. in 
$X_{t}$ of length $k$
 with constant $1+\epsilon_0$.

To prove this Lemma we shall apply Lemma 7. First,
 we show that if $G_1^*,\ldots,G_k^*$
 is any special sequence
 in $\Delta$ of length $k$ and
$E$ is any interval, then $|z^*(Ex)|<1$,
 where $z^*$ is the $(k,\sqrt
f)$-form $f(k)^{-1/2}\sum_{i=1}^k z_i^*$
with $z_j^*=G_j^*(it)$, and $x=\sum_{i=1}^k x_i$.

Indeed, let
 $s$ be maximal such that $G_s^*=F_s^*$ or zero
 if no such $s$
exists.
Suppose now $i\neq j$ or one of $i,j$ is greater than $s+1$.
 Then since
$\sigma$ is an injection,
 we can find $L_1\neq L_2\in L$ such that $z_i^*$ is an
$(L_1,f)$-form and $x_j$ is the normalized sum of a R.I.S.
 of length $L_2$ and
also an $l_{1+}^{L^{\prime}_2}$-average with
 constant $1+\epsilon_0$, where
$L^{\prime}_2=L_2^{\epsilon_0 /4}$.
 We can now use Lemmas~GM4 and~GM5 to show that
 $|z_i^*(E x_j)|<k^{-2}$.

 If $L_1 < L_2$,
 it follows from the
lacunarity of $L$ that $L_1 < L^{\prime}_2$.
 We know that $L_1 \geq j_{2k}$ since
$L_1$ appears in a special sequence of length $k$.
 Lemma~GM4 thus gives
 $|z_i^*(Ex_j)| = |(Ez^*_i)(x_j)| \leq 3(1+\epsilon_0)/f(L_1)$.
 The conclusion in this
case now follows from the fact that
 $f(l)\geq 4 k^2$ when $l\geq j_{2k}$.

If $L_2 < L_1$, we apply Lemma~GM5 in $X_{t}$ with $\epsilon = 1$
 to the non-normalized
sum $x^{\prime}_j$ of the R.I.S. the normalized sum of which is $x_j$.
 The definition of
$L$ gives us that $M_f(L_2) < L_1$, so Lemma~GM5 gives
 $|z^*_i(Ex^{\prime}_j)| \leq 3$.
It follows from the lower $f$-estimate in $X_{t}$ that
 $\|x^{\prime}_j\| \geq L_2/f(L_2)$. The
conclusion now follows because $l\geq j_{2k}$
 implies that $f(l)/l \leq 1/4k^2$.

\

Now choose an interval $E^{\prime}$ such that
\[ \Bigl|(\sum_{i=1}^s z^*_i)(Ex)\Bigl| =
 \Bigl|(\sum_{i=1}^k F^*_i(it))(E^{\prime} x)\Bigl| \leq 4.\]
It follows that
\[ \Bigl|(\sum_{i=1}^k z^*_i)(Ex)\Bigl| \leq 4 +|z_{s+1}^*(x_{s+1})| +
k^2.k^{-2} \leq  6 .\]
We finally obtain that $|z^*(Ex)|\leq 6f(k)^{-1/2} < 1$ as claimed.

Now let $\phi^{\prime}$ be the function
\[ \phi^{\prime}(x)=\sqrt{f(x)} \ if \  x\in K, x \neq k,
   \phi^{\prime}(x)=f(x) \ otherwise. \]
Let $g^{\prime}$ be the function obtained from $\phi^{\prime}$ 
by Lemma 8 in the case
$K_{0}=K \setminus \{k\}$; we know that $g^{\prime}(l)=f(l)$ for
 every $l \in L \cup \{k\}$.

It follows from what we have just shown about special sequences of length $k$
that for every interval $E$,
\[ \|Ex\|_{t} \leq 1\vee\sup\{|x^*(Ex)|:
M\geq 2, x^*\ is\ an\ (M,g^{\prime})-form\}
.\]
Since $x$ is the sum of a
R.I.S., Lemma~7 implies that $\|x\|_{t} \leq
(1+2\epsilon_0)k g^{\prime}(k)^{-1} = (1+2\epsilon_0) k/f(k)$.

\section{ $X$ is hereditarily indecomposable }

Let $Y$ and $Z$ be two infinite-dimensional subspaces of $X$.
We want to show that their sum is not a topological sum.
 Let \-$\delta>0$. We shall build two vectors $y \in Y$ and $z \in Z$
such that $\delta \|y+z\| > \|y-z\|$.

Let $\epsilon_{0}= 1/10$.
 Let $k \in K $
 be an integer such that
 $1/4 <\epsilon_{0}\ k^{1/p}/f(k)^{1-\theta}$ and
 $2/\sqrt{f(k)}^{1-\theta} < \delta$, 
and let $\epsilon >0 $ be such that  $\sqrt{\epsilon} \leq \epsilon_{0} / 4kf(k) $.
We may assume that both $Y$ and $Z$ are spanned by block bases.
By Lemma~2,
$Y$ and $Z$ contain, for every $N \in \N$,
 an $l^{N}_{p+}$-average with constant $1+\epsilon$. We now
build a sequence $(x_{j})_{j=1}^{k}$ in $X$ by iteration.

\paragraph{First step}
Let $x_{1} \in Y$ be a R.I.S.-vector of norm 1, constant $1+\epsilon$ and length
$M_{1}=j_{2k}$; we have $M_{1}^{\epsilon_{0}/4} = N_{1} \geq 4M_{f} (k/\epsilon_{0})/
\epsilon_{0} \ff(1) $.
Let $x_{11} < \cdots < x_{1M_{1}}$ be the R.I.S.
 whose normalized sum is $x_{1}$: there exists $\lambda_{1}$ such that
$ \lambda_{1} x_{1}=x_{11}+\cdots+x_{1M_{1}}$.
Applying the lower estimate in $X$
and Lemma 10, we get
\[ M_{1}^{1/p} / f(M_{1})^{1-\theta} \leq
\| \lambda_{1} x_{1} \| \leq
(1 + 10 \sqrt{\epsilon}) M_{1}^{1/p} / f(M_{1})^{1-\theta}.\]
so that $\lambda_{1} = M_{1}^{1/p} / 
f(M_{1})^{1-\theta} $ up to the multiplicative factor
 $1 + 10 \sqrt{\epsilon}$.

Now we associate to $x_{1m}$:
\begin{itemize}
\item a representative function $F_{1m}$ for $x_{1m}$;
\item a vector $x_{1m}^{*}$ in $X^{*}$ that norms $x_{1m}$ and with
 $ran(x_{1m}^{*}) \subset ran(x_{1m})$;
\item a minimal interpolation function
 $F_{1m}^{*}$ for $x_{1m}^{*}$; it exists because
 of Proposition 2 and because, as 
$x_{1m}^{*}$ is of finite range, Theorem 1 applies.
\end{itemize}

The function $F_{1m}^{*}$ is in $\overline{\cal S}$. Indeed, remember
that if we
restrict them to finite range vectors, the unit
ball of $X_{t}^{*}$ and $\overline{D(t)}$ coincide; so by the
convexity of $D(t)$, for
every $\nu >0$, the function $F_{1m}^{*}/(1+\nu)$ takes its values
in $D(t)$ for almost every $it$; as it takes its values in
 $B(l_{\qiq}) \ a.e.$ on $S_{1}$, it is in $\cal S$, which ends
the proof.

 By Lemma 5, $F_{1m}^{*}$ can be approached by a
 function ${\cal F}_{1m}^{*}$
in $\Delta$  ( and because of the interval projection stability of
$\Delta$, we may assume that
 $ran({\cal F}_{1m}^{*}) \subset ran(F_{1m}^{*})$).
 More precisely, we suppose that
${\cal F}_{1m}^{*}$ is close to $F_{1m}^{*}$
 up to $\epsilon / (1+\epsilon) $
 for the norm $\int_{z \in \delta S} \|.\|^{*}_{z} d\mu(z)$
 ( over functions of finite range, this norm is equivalent
to the norm $\int_{z \in \delta S} \|.\|_{1} d\mu(z)$
 first introduced).

Lastly, we define two functions:

Let ${\cal F}_{1}^{*} = f(M_{1},z)^{-1}
 \sum_{m=1}^{M_{1}} {\cal F}_{1m}^{*}$. It belongs to $\Delta$.
Let $x_{1}^{*}= {\cal F}_{1}^{*} (\theta)$.

Let  $F_{1} = f(M_{1},z) M_{1}^{-1}
 \sum_{m=1}^{M_{1}} F_{1m}$.

\paragraph{Iteration}
Let $M_{2}=\sigma({\cal F}_{1}^{*}) \in L$.
We may assume we chose ${\cal F}_{1}^{*}$ such that
 $ 1/2 \ \epsilon_{0} f(M_{2}^{\epsilon_{0}/4})^{1/2}
 \geq |ran(F_{1})| $ ( by choosing a function ${\cal F}^{*}_{11}$ such
that $M_{2}$ is big enough).
 Let $x_{2} \in Z$ be a R.I.S.-vector of norm $1$,
 constant $1+\epsilon$ and length
$M_{2}$, and $x_{2} > x_{1}$; and repeat the above construction. By iterating it, we 
obtain for
$j=1,\ldots,k$ sequences $F_{j},x_{j},{\cal F}_{j}^{*},x_{j}^{*}$ such that:

\begin{itemize}

\item $x_{j} \in Y$ when $j$ is odd, $x_{j} \in Z$ otherwise.
\item $\|x_{j}\|=1 $ for every $j$ and $\|x_{j}^{*}\|^* \leq 1 $.
\item $x_{j}=F_{j}(\theta)$ up to $10 \sqrt{\epsilon}$ and
 $x_{j}^{*}={\cal F}_{j}^{*}(\theta)$.
\item ${\cal F}_{1}^{*},\ldots,{\cal F}_{k}^{*}$ is a
 special sequence of length $k$.
\item For $j=2,\ldots,k$,
 $1/2 \ \epsilon_{0} f(M_{j}^{\epsilon_{0}/4})^{1/2}
 \geq |ran(F_{j-1})|. $
\item For every $j$, $< {\cal F}_{j}^{*}(\theta),F_{j}(\theta) > = 1 $ up to 
$\epsilon$.
\item For every $j$, except on ${\cal J}_{j}$ of measure at
 most $2 \sqrt{\epsilon}$,
 $< {\cal F}_{j}^{*}(it),F_{j}(it) >=1 $ up to
 $2 \sqrt{\epsilon}$.
\item For every $j$, except on ${\cal J}^{\prime}_{j}$
 of measure at most $4 \sqrt{\epsilon} / f(M_{j})$, $F_{j}(it)$ is up to
 $10 \sqrt{\epsilon}$ the normalized sum
of a R.I.S. with constant $1+4\sqrt{\epsilon} \leq 1+\epsilon_{0}/4$.
\end{itemize}

{\em Proof } Only the last three points are not obvious.

\

\noindent{\em First point }
For $F$ and $F^{*}$ in ${\cal A}_S$, define 
$<F,F^{*}>$ to be $\int_{z \in \delta S} <F(z),F^{*}(z)> d\mu(z)$, and
notice that this is equal to $<F(\theta),F^{*}(\theta)>$ by
 analyticity. Now for every $j$,
\[ < {\cal F}_{j}^{*},F_{j}> =
 \frac{1}{M_{j}} \sum_{m=1}^{M_{j}} < {\cal F}_{jm}^{*},F_{jm}>. \]
If we replace each ${\cal F}_{jm}^{*}$ by $F_{jm}^{*}$, the sum
is equal to \[\frac{1}{M_{j}} \sum_{m=1}^{M_{j}} <x^{*}_{jm},x_{jm}> = 1.\]

The error we make by doing this is
$| 1/M_{j} \sum_{m=1}^{M_{j}} 
<{\cal F}_{jm}^{*}-F_{jm}^{*},F_{jm}>|$, smaller than
\[ \frac{1}{M_{j}} \sum_{m=1}^{M_{j}} 
\int_{z \in \delta S}
\| ({\cal F}_{jm}^{*}-F_{jm}^{*})(z) \|_{z}^{*}
\| F_{jm}(z)\|_{z} d\mu(z) 
\leq \frac{1}{M_{j}} \sum_{m=1}^{M_{j}} \frac{\epsilon}{1+\epsilon}
 (1+\epsilon)
\leq \epsilon \]
(we recall that as $F_{jm}$ is representative
 for $x_{jm}$, $\|F_{jm}(z)\|_{z} \leq 1+\epsilon \ a.e.$).

\

{\em Second point } Let $F_{j}^{*}$ be the function 
$f(M_1,z)^{-1}
 \sum_{m=1}^{M_{1}} F_{jm}^{*}$. It is easy to see that

\[ 1=\int_{z \in \delta S} <F_{j}^{*}(z),F_{j}(z)> d\mu(z),\]
while
\[<F_{j}^{*}(z),F_{j}(z)> \leq 1+\epsilon\ a.e.\ . \]

By a Bienaym\'{e}-Tchebitschev estimation, except on a set of measure
at most $\sqrt{\epsilon}$, $<F_{j}^{*}(z),F_{j}(z)>=1$ up to
$\sqrt{\epsilon}$. Furthermore, we know that

\[ \int_{z \in \delta S}
 |<({\cal F}_{j}^{*}-F_{j}^{*})(z),F_{j}(z)>| d\mu(z)  \leq \epsilon\]

so that except on a set of measure at most $\sqrt{\epsilon}$,
$<({\cal F}_{j}^{*}-F_{j}^{*})(z),F_{j}(z)>=0$ up
 to $\sqrt{\epsilon}$.

Adding these two estimates completes the proof.

\

{\em Third point }
For each $m$, $F_{jm}$ is representative for $x_{jm}$, so by Lemma 4,
except on a set ${\cal J}^{\prime}_{j}$ of measure
 $4\sqrt{\epsilon}/f(M_{j})$, we have
 \[ \Bigl\| \sum_{m=1}^{M_{j}}F_{jm}(it)-x_{t} \Bigl\|_{t}
 \leq 2\sqrt{\epsilon}, \]
where $x_{t}$ is the sum of a R.I.S. in $X_{t}$ with constant
 $1+4\sqrt{\epsilon}$. So
\[ \left\| F_{j}(it) - \left(\frac{M_{j}^{1/\qiq}}{f(M_{j})}\right)^{it} 
\frac{f(M_{j})}{M_{j}}x_{t} \right\|_{t} \leq
2\sqrt{\epsilon} \frac{f(M_{j})}{M_{j}} \leq
2\sqrt{\epsilon}. \]
The proof follows, because
 by Lemma 9, 
$$ M_{j}/f(M_{j})\leq \|x_{t}\|_{t}
\leq (1+8\sqrt{\epsilon})M_{j}/f(M_{j}),$$
 so $f(M_{j})/M_{j}\ x_{t}$ is up to $8\sqrt{\epsilon}$ a
{\em normalized} R.I.S.-vector.

\paragraph{Estimation of $\|\sum_{j=1}^{k} x_{j}\|$}

Let ${\cal G}^{*} =  g(k,z)^{-1}  
\Sigma_{j=1}^{k}{\cal F}_{j}^{*}$. Since for every $j$,
${\cal F}_{j}^{*} \in \Delta$, and $k$ is in $K$,  ${\cal G}^{*}$
 is in $\cal D$ and by Lemma 6, $x^{*}={\cal G}^{*}(\theta)$ is in
 the unit ball of $X^{*}$.

So $\| \sum_{j=1}^{k} F_{j}(\theta) \| 
\geq x^{*}( \sum_{j=1}^{k} F_{j}(\theta)) \geq (1-\epsilon) k^{1/p}/ 
\sqrt{f(k)}^{1-\theta} $, and
\[\Bigl\| \sum_{j=1}^{k} x_{j} \Bigl\| \geq (1-\epsilon_{0}) k^{1/p}/ 
\sqrt{f(k)}^{1-\theta}- 1/4 \geq
(1-2\epsilon_{0}) k^{1/p}/ 
\sqrt{f(k)}^{1-\theta}. \]
( the 1/4 is the error we made by replacing the $x_{j}$'s
 by the $F_{j}(\theta)$'s ).

\paragraph{Estimation of $\|\sum_{j=1}^{k} (-1)^{j-1} x_{j}\|$}

Let ${\cal J}$ be the union of the ${\cal J}_{j}$'s and the 
${\cal J}^{\prime}_{j}$'s.
 The set $\cal J$ is of measure at most
$6k\sqrt{\epsilon}$.

For every $t$ in ${\R} \setminus {\cal J}$, for
 every interval $E$, let us evaluate

\[ \left|\left(\sum_{j=1}^{k} {\cal F}_{j}^{*}(it)\right)
\left(\sum_{j=1}^{k} (-1)^{j-1} EF_{j}(it)\right)\right|. \]

This is a sum of at most $k$ scalars. Those who come
 from terms of range included in
 $E$
are equal to $(-1)^{j-1}$ up to $2\sqrt{\epsilon}$, so that their sum
is $-1,0$ or $1$ up to $2k\sqrt{\epsilon}$; two others can come from
terms whose range intersects $E$, they are bounded in modulus
 by $1+10\sqrt{\epsilon}$; the others are equal to $0$. So the sum
is smaller than $1+2k\sqrt{\epsilon}+2(1+10\sqrt{\epsilon})
\leq 3+3k\sqrt{\epsilon}$.

For every $j$, $F_{j}(it)$ is up to $10\sqrt{\epsilon}$ a
 R.I.S. vector $x_{j}(t)$. The $(-1)^{j-1}x_{j}(t)$'s satisify
 the hypotheses of Lemma 12:
the increasing condition is satisfied, and 
for every interval $E$,
\[  \left|\left(\sum_{j=1}^{k} {\cal F}_{j}^{*}(it)\right)
\left(\sum_{j=1}^{k}(-1)^{j-1}Ex_{j}(t)\right)\right|
\leq 3+3k\sqrt{\epsilon}+10k\sqrt{\epsilon} \leq 4. \]
 It then follows from the conclusion of Lemma 12 and the relation
between $F_{j}(t)$ and $x_{j}(t)$ that  
\[ \Bigl\|\sum_{j=1}^{k} (-1)^{j-1} F_{j}(it) \Bigl\|_{t} \leq
 (1+2\epsilon_{0})k/f(k)+10k\sqrt
{\epsilon} .\]
It follows that
\[\int_{\R \setminus {\cal J}}
 \Bigl\|\sum_{j=1}^{k} (-1)^{j-1} F_{j}(it) \Bigl\|_{t}
 d\mu_{0}(t) \leq
 (1+2\epsilon_{0})k/f(k)+10k\sqrt
{\epsilon} .\]

\

We now want to estimate the integral of this same norm on ${\cal J}$. 
It is enough, by a triangular inequality, to evaluate
$\int_{t \in {\cal J}} \| F_{j}(it) \|_{t}
 d\mu_{0}(t)$.
 If $t$ belongs to ${\cal J}^{\prime}_{j}$, by
a triangular inequality, 
$\| F_{j}(it) \|_{t} \leq (1+\epsilon) f(M_{j})$, but recall
that ${\cal J}^{\prime}_{j}$ is of measure at most 
$4\sqrt{\epsilon}/ f(M_{j})$; else, $F_{j}(it)$ is up to
$10\sqrt{\epsilon}$ a
normalized R.I.S. vector,
 so that $\| F_{j}(it) \|_{t} \leq 1+10\sqrt{\epsilon}$, and this
on a set of measure less than $6k\sqrt{\epsilon}$. Finally,
\[ \int_{\cal J} \| F_{j}(it) \|_{t}
 d\mu_{0}(t) \leq
 6k\sqrt{\epsilon}(1+10\sqrt{\epsilon})+
\frac{4 \sqrt{\epsilon}}{f(M_{j})}(1+\epsilon)f(M_{j})
\leq 7k\sqrt{\epsilon}. \]
and 
\[\int_{\cal J}
 \Bigl\|\sum_{j=1}^{k} (-1)^{j-1} F_{j}(it) \Bigl\|_{t}
 d\mu_{0}(t) \leq
 7k^{2}\sqrt{\epsilon}.\]

\

It follows from these two estimates that
\[\int_{\R}
 \Bigl\|\sum_{j=1}^{k} (-1)^{j-1} F_{j}(it) \Bigl\|_{t}
 d\mu_{0}(t) \leq
 (7k^{2}+10k)\sqrt{\epsilon}+(1+2\epsilon_{0})\frac{k}{f(k)}
\leq (1+4\epsilon_{0}) \frac{k}{f(k)}.\]
Furthermore, almost everywhere on $S_{1}$,
 \[ \Bigl\|\sum_{j=1}^{k} (-1)^{j-1} F_{j}(1+it) \Bigl\|_{q}
 \leq (1+\epsilon)k^{1/q} ,\]
so that, by Lemma 1,

\[\Bigl\| \sum_{j=1}^{k} (-1)^{j-1} F_{j}(\theta) \Bigl\|
 \leq (1+3\epsilon_{0}) k^{1/p}/ 
f(k)^{1-\theta} ,\]
and
\[\Bigl\| \sum_{j=1}^{k} (-1)^{j-1} x_{j} \Bigl\|
 \leq (1+3\epsilon_{0}) k^{1/p}/ 
f(k)^{1-\theta}+1/4 \leq
(1+4\epsilon_{0}) k^{1/p}/ 
f(k)^{1-\theta}.\]

\paragraph{Conclusion}
Let $y \in Y$ be the sum of the $x_{j}$ with odd indices, $z \in Z$ be the sum 
of the $x_{j}$ with even indices. By the above estimates and by
choice of $k$, they satisfy $\delta \|y+z\| > \|y-z\|$.
 As $\delta$ is arbitrary and so are $Y$ and $Z$, $X$ is 
hereditarily indecomposable.

\paragraph{}I warmly thank Bernard Maurey for his help.

\end{document}